\documentstyle[amsfonts,12pt,thmsa,sw20lart]{article}

\input tcilatex
\begin{document}

\title{{\bf A merry-go-round with the circle map, primes and pseudoprimes }$^{1}$}
\author{M. Leo, R.A. Leo, and G. Soliani \\
Dipartimento di Fisica dell'Universit\`{a} di Lecce, 73100 Lecce, Italy, \\
and Istituto Nazionale di Fisica Nucleare, Sezione di Lecce, Italy}
\date{}
\maketitle

{\bf Abstract}. We show that the use of the main characteristics of the
circle map leads naturally to establish a few statements on primes and
pseudoprimes. In this way a Fermat's theorem on primes and some interesting
properties of pseudoprimes are obtained.

\section{Introduction}

The theory of prime numbers (primes) attracted (and continue to attract)
both mathematicians and physicists, and computer scientists$^{2,3}$.
Although a lot of works dwell upon problems concerning this topic, a minor
attention is payed on the possible connections between the theory of primes
and other subjects. For instance, factorization techniques into primes turn
out to be important in the construction of codes for transmitting reserved
information$^{3}$. Primes appear also in the investigation of quantum chaos
of dynamical systems$^{4}$. In fact, in this context analogies were found
between classical periods and primes, and energy levels and zeros of
Riemann's zeta function. On the other hand, concerning the application of
primes to physics, a relation between prime sequences and the spectra of
excited nuclei have been suggested$^{4}$.

Hence, the search for the existence of links involving primes and other
branches of physics and mathematics could be a challenge for getting further
insights into the theory of primes.

Following this spirit, here we show that some properties of primes and
psudoprimes can be derived in a simple way starting from the circle map. By
virtue of its pedagogical valence, this approach could be useful mostly to
beginners.

In Sec. 2 some preliminaries on the circle map are expounded. Furthermore,
by using the properties of the periodic orbits of this map a well-known
Fermat's theorem on primes is formulated. In Sections 3 and 4 we show as the
characteristics of the circle map can be usefully employed also to obtain
some properties of pseudoprimes. Finally, in Sec. 5 some concluding remarks
are reported.

\section{The circle map and a Fermat's theorem on primes}

If $S^{1}$ is the unitary circonference, let us define the application $f:$ $%
S^{1}\longrightarrow S^{1}$ as

\begin{equation}
f(\theta )=k\theta \text{ }(\func{mod}2\pi ),  \label{1}
\end{equation}

where $\theta \in [0,2\pi ]$ denotes the angular position of a point
belonging to $S^{1}$ and $k\in {\Bbb N}:$ $k>1.$ The application $f$ is
called circle map$^{5}$; it establishes a correspondence between points
lying on $S^{1}.$

The fixed points of the map (1) are obtained by solving the equation

\begin{equation}
f(\theta )=\theta \text{ }(\func{mod}2\pi )  \label{2}
\end{equation}

which gives

\begin{equation}
\theta _{j}=\frac{2j}{k-1}\pi ,  \label{3}
\end{equation}

with $j=0,1,...,k-2.$ There are $k-1$ fixed points.

On the other hand, the periodic points of (1) of period $n\in {\Bbb N}$ can
be determined from the equation

\begin{equation}
f^{n}(\theta )=\theta \text{ }(\func{mod}2\pi )  \label{4}
\end{equation}

which yields

\begin{equation}
\theta _{j}=\frac{2j}{k^{n}-1}\pi ,  \label{5}
\end{equation}

with $j=0,1,...,k^{n}-2.$

We notice that the fixed points are also periodic points of period $n\in 
{\Bbb N}:n>1.$ Then, excluding the fixed points, the number of periodic
points of period $n$ is

\begin{equation}
k^{n}-1-(k-1)=k^{n}-k.  \label{6}
\end{equation}
Now, if $n$ is a prime, the periodic points of period $n$ (excluding the
fixed points), are {\it necessarily} periodic points of {\it prime }period $%
n.$

In other words, any periodic orbit of {\it prime} period{\it \ }$n$ is
constituted by $n$ points. Then

\begin{equation}
\frac{k^{n}-k}{n}  \label{7}
\end{equation}
($n$ being a prime) is equal to the number of periodic orbits of (prime{\it %
) }period $n$ and, consequently, is an integer.

Thus, we can formulate Fermat's theorem on primes$^{2,3}$ as follows:

\smallskip

$\bullet $ {\it If }$n\in {\Bbb N}$ {\it is a prime} ($n>1)${\it \ and }$%
k\in {\Bbb N}:k>1,${\it \ then }

\begin{equation}
\frac{k^{n}-k}{n}\in {\Bbb N}.  \label{8}
\end{equation}
In particular:

$\bullet $ {\it If }$k${\it \ and }$n${\it \ are relatively primes, i.e. GCD
(k,n) =1, (8) is equivalent to}

\begin{equation}
\frac{k^{n-1}-1}{n}\in {\Bbb N}.
\end{equation}

\section{The circle map and pseudoprimes}

\smallskip We remind the reader the definition of pseudoprime (see, for
example, Ref. 3, p. 33):

$\bullet $ If $n$ is an odd composite number which is relatively prime to $k$
and if

\smallskip

\begin{eqnarray}
k^{n-1} &\equiv &1\text{ }(\func{mod}n),  \label{10} \\
&&  \nonumber
\end{eqnarray}

then $n$ is called a {\it pseudoprime of basis }$k.$ In other words, (10)
tells us that

\begin{equation}
\frac{k^{n}-k}{n}\in {\Bbb N}.
\end{equation}

In the definition (10) the notation $a\equiv b$ $(\func{mod}n)$ means that $%
n $ divides $a-b$, i.e. $\frac{a-b}{n}\in {\Bbb N}$ (The smallest
pseudoprime of basis 2 is 341=11$\times 31).$

We recall that there exist composite numbers which are pseudoprimes for all
bases to which they are relatively prime, in the sense specified in Sec. 4.
The first of such numbers, called {\it Carmichael numbers}, is $561=3\cdot
11\cdot 17$ (see Ref. 3, p. 33).

\smallskip We have the following

{\bf Theorem 1} {\it If }$n${\it \ is a pseudoprime of basis }$k,${\it \ then%
}

\begin{equation}
\frac{k^{n}-k-\Pi _{n}}{n}\in {\Bbb N},  \label{11}
\end{equation}
{\it where }$\Pi _{n}${\it \ stands for the number of periodic points of
prime period }$n${\it \ }$(n\in {\Bbb N})${\it \ in the circle map (1).}

The property (12) emerges immediately observing that $\Pi _{n}$ is always
divisible by $n.$

\smallskip

\section{Factorization of pseudoprimes}

\smallskip In this Section we shall prove some interesting properties of
pseudoprimes which turn out to be expressed as products of two and three
primes, respectively. We recall that the decomposition of a positive integer
into primes is named its {\it factorization}. Namely, this can be written as

\smallskip

\begin{equation}
n=n_{1}^{a_{1}}\times n_{2}^{a_{2}}\times \cdot \cdot \cdot \times
n_{r}^{a_{r}},
\end{equation}

where $n_{1},n_{2},...,n_{r}$ are distinct primes and $a_{1},a_{2},...,a_{r}$
$\in {\Bbb N}.$

Moreover, if the factorizations of two integers $n$ and $m$ have no common
primes, then $n$ and $m$ are said {\it relatively }prime.

\subsection{Pseudoprimes product of two primes}

Let us suppose now that $n$ be a pseudoprime of basis $k$ of the type

\begin{equation}
n=n_{1}n_{2},  \label{12}
\end{equation}
\[
\]

where $n_{1}$ and $n_{2}$ are distinct primes. We have

\smallskip

\begin{equation}
k^{n}-k-\Pi _{n}=\Pi _{n_{1}}+\Pi _{n_{2}},
\end{equation}

\smallskip

where $\Pi _{n}$ is the number of periodic points of prime period $n,$ and

\begin{equation}
\Pi _{n_{1}}=k^{n_{1}}-k,  \label{13}
\end{equation}

\begin{equation}
\Pi _{n_{2}}=k^{n_{2}}-k,  \label{14}
\end{equation}

denote the numbers of periodic points of prime periods $n_{1}$ and $n_{2},$
respectively.

The fact that $n$ is a pseudoprime of basis $k$ implies the validity of
(12). In other words, the sum of the number of the periodic points of prime
period $n_{1}$ and the number of the periodic points of prime period $n_{2}$
is divisible by $n.$ Indeed

\smallskip

\begin{equation}
\frac{\Pi _{n_{1}}+\Pi _{n_{2}}}{n}\in {\Bbb N}\Longleftrightarrow \frac{%
k_{{}}^{n_{1}}-k+k^{n_{2}}-k}{n}\in {\Bbb N}\text{ }  \label{16}
\end{equation}

or, equivalently:

\[
\]

\begin{equation}
\frac{k_{{}}^{n_{1}}-k+k^{n_{2}}-k}{n_{1}}\in {\Bbb N},\text{ }\frac{%
k_{{}}^{n_{1}}-k+k^{n_{2}}-k}{n_{2}}\in {\Bbb N}\text{.}  \label{17}
\end{equation}
Furthermore, since $n_{1}$ and $n_{2}$ are primes, namely $\frac{k^{n_{1}}-k%
}{n_{1}}\in {\Bbb N},$ $\frac{k^{n_{2}}-k}{n_{2}}\in {\Bbb N},$ then

\smallskip

\begin{equation}
\frac{k^{n_{2}}-k}{n_{1}}\in {\Bbb N},\text{ }\frac{k^{n_{1}}-k}{n_{2}}\in 
{\Bbb N}.  \label{18}
\end{equation}
Now we are ready to enounce the following

\smallskip

{\bf Theorem 2} {\it Let }$n_{1},n_{2}$ {\it be\ two distinct primes and }$%
k\in {\Bbb N}:k>1${\it \ a number such that }$n_{1},n_{2}${\it \ and }$k$%
{\it \ are relatively primes, in the sense that }$GCD$ $(n_{1},k)=GCD$ $%
(n_{2},k)=1.${\it \ Then, }$n=n_{1}n_{2}${\it \ is a pseudoprime of basis }$%
k ${\it \ if and only if the conditions (20) hold.} 
\[
\]

{\bf Proof}

Let us assume that $n_{1},n_{2}$ and $k$ are relatively primes. We have just
shown that if $n=n_{1}n_{2}$ (where $n_{1},n_{2}$ are two distinct primes)
is a pseudoprime of basis $k,$ then the properties (20) are valid.
Vice-versa, if the relations (20) are satisfied, then, since by hypothesis $%
n_{1},n_{2}$ are two primes, we have

\smallskip

\begin{equation}
\frac{k^{n_{1}}-k}{n_{2}}\in {\Bbb N}\Longrightarrow \frac{k^{n_{1}}-k}{%
n_{1}n_{2}}\in {\Bbb N},  \label{20}
\end{equation}

\begin{equation}
\frac{k^{n_{2}}-k}{n_{1}}\in {\Bbb N}\Longrightarrow \frac{k^{n_{2}}-k}{%
n_{1}n_{2}}\in {\Bbb N}.  \label{21}
\end{equation}
The conditions (21) and (22) imply

\smallskip

\begin{equation}
\frac{k^{n_{1}}-k+k^{n_{2}}-k}{n_{1}n_{2}}\in {\Bbb N}.  \label{22}
\end{equation}
Now, keeping in mind the relations (15), (16) and (17) we deduce that $\frac{%
k^{n}-k}{n}\in {\Bbb N},$ namely $n$ is a pseudoprime of basis $k.$

{\bf Remark}

We have seen that if $n=n_{1}n_{2},$ with $n_{1},n_{2}$ primes, is a
pseudoprime of basis $k,$ then

\smallskip

\begin{equation}
\frac{k^{n_{1}}-k}{n_{{}}}\in {\Bbb N},\text{ }\frac{k^{n_{2}}-k}{n_{{}}}\in 
{\Bbb N}.  \label{27}
\end{equation}

Now, if $n_{1}>n_{2},$ the relations (24) provide

\smallskip

\begin{equation}
\frac{k^{n_{1}}-k^{n_{2}}}{n_{{}}}\in {\Bbb N},  \label{28}
\end{equation}
from which

\smallskip

\begin{equation}
\frac{k^{n_{1}-n_{2}}-1}{n_{{}}}\in {\Bbb N}.  \label{29}
\end{equation}

\smallskip

This condition is equivalent to the two relations:

\smallskip

\begin{equation}
\frac{k^{n_{1}-n_{2}}-1}{n_{1}}\in {\Bbb N},\text{ }\frac{k^{n_{1}-n_{2}}-1}{%
n_{2}}\in {\Bbb N}.  \label{30}
\end{equation}

\subsubsection{Generalizations of the properties (27)}

\smallskip

The properties (27) can be generalized as follows.

{\bf a)} Let us note that

\smallskip

\begin{equation}
k^{n_{1}^{r}}-k=\stackrel{r}{\stackunder{j=1}{\sum }}%
(k^{n_{1}^{j}}-k^{n_{1}^{j-1}}),  \label{31}
\end{equation}
$\forall r\in {\Bbb N}$.

Thus, since

\begin{equation}
\frac{k^{n_{1}^{j}}-k^{n_{1}^{j-1}}}{n_{1}}\in {\Bbb N},\text{ }j=1,2,...,r,
\label{32}
\end{equation}
we get

\begin{equation}
\frac{k^{n_{1}^{r}}-k}{n_{1}}\in {\Bbb N}.  \label{33}
\end{equation}
Combining together the relation (30) and the second of (24), we find

\smallskip

\begin{equation}
\frac{k^{n_{1}^{r}}-k^{n_{2}}}{n_{1}}\in {\Bbb N}\Longleftrightarrow \frac{%
k^{n_{1}^{r}-n_{2}}-1}{n_{1}}\in {\Bbb N}.  \label{34}
\end{equation}

\smallskip In a similar way we can prove the property

\smallskip

\begin{equation}
\frac{k^{\mid n_{2}^{r}-n_{1}\mid }-1}{n_{2}}\in {\Bbb N},\text{ }\forall
r\in {\Bbb N}.  \label{35}
\end{equation}

\[
\]
We observe that the exponent of $k$ has to be positive; this justifies the
appearance of the absolute value in (32).

\smallskip

${\bf b)}$ Let $n=n_{1}n_{2}$ (where $n_{1},n_{2}$ are primes) be a
pseudoprime (of basis $k$ ) relatively prime to $n_{1}$ and $n_{2}.$ We have
shown that (see (24))

\smallskip

\begin{equation}
\frac{k^{n_{i}-1}-1}{n}\in {\Bbb N},\text{ }i=1,2.  \label{36}
\end{equation}
The relation (31) entails

\smallskip

\begin{equation}
k^{n_{i}-1}\frac{k^{n_{i}-1}-1}{n}\in {\Bbb N}\Longrightarrow \frac{%
k^{2(n_{i}-1)}-k^{n_{i}-1}}{n}\in {\Bbb N}\Longrightarrow \frac{%
k^{2(n_{i}-1)}-1}{n}\in {\Bbb N}.  \label{37}
\end{equation}
By iterating the procedure, we find 
\begin{eqnarray*}
&& \\
&&
\end{eqnarray*}
\begin{equation}
\frac{k^{r(n_{i}-1)}-1}{n}\in {\Bbb N},\text{ }\forall r\in {\Bbb N}(i=1,2).
\label{38}
\end{equation}

$\smallskip $

$\bullet $ From (31) we can derive an interesting result, which can be
deduced also from Euler's Theorem (see, for example, Ref. 3, p. 35).
Precisely, (31) gives

\smallskip 
\begin{equation}
k^{n_{2}-1_{{}}}\frac{k^{n_{1}-1}-1}{n}\in {\Bbb N}\Longrightarrow \frac{%
k^{n_{1}+n_{2}-2}-k^{n_{2}-1}}{n}\in {\Bbb N}\Longrightarrow \frac{%
k^{n_{1}+n_{2}-2}-1}{n}\in {\Bbb N}.  \label{39}
\end{equation}

\[
\]
This result can be obtained simply exploiting Euler's Theorem, which we
report below for reader's convenience:

\smallskip

{\bf Euler's Theorem}

\smallskip

{\it Let }$\phi (n)${\it \ denote the number of positive integers less than
or equal to }$n${\it \ and relatively prime to }$n.${\it \ Then, if }$n${\it %
\ and }$k${\it \ are positive and relatively prime integers, the property}

\begin{equation}
k^{\phi (n)}\equiv 1\text{ }(\func{mod}n)  \label{40}
\end{equation}

{\it holds.}

We notice that if $n$ is prime, then $\phi (n)=n-1.$

\smallskip

To consider the relation (34) as a special case of Euler's Theorem, it is
sufficient to keeping in mind that when $n=n_{1}n_{2}$ (with $n_{1}$ and $%
n_{2}$ primes), then $\phi (n)=n_{1}n_{2}$ $-n_{1}-n_{2}+1.$

\smallskip

In fact, since $\frac{k^{^{n_{1}n_{2}-1}}-1}{n}\in {\Bbb N},$ we have

\smallskip

\begin{equation}
\frac{k^{n_{1}n_{2}-n_{1}-n_{2}+1.}-1}{n}\in {\Bbb N}\Longleftrightarrow 
\frac{k^{n_{1}n_{2}-1}-k^{n_{1}+n_{2}-2}}{n}\in {\Bbb N}\Longleftrightarrow 
\frac{k^{n_{1}+n_{2}-2.}-1}{n}\in {\Bbb N}.  \label{41}
\end{equation}
${\bf c)}$ An extension of (36) is represented by

\smallskip

\begin{equation}
\frac{k^{rn_{1}+sn_{2}-(r+s).}-1}{n}\in {\Bbb N},\text{ }\forall r,s\in 
{\Bbb N}.  \label{42}
\end{equation}
Indeed, with the help of (35) we can write, for $\forall r\in {\Bbb N}:$

\smallskip

\begin{equation}
\frac{k^{rn_{1}-r}-1}{n}\in {\Bbb N}\Longrightarrow k^{n_{2}-1}\frac{%
k^{rn_{1}-r}-1}{n}\in {\Bbb N}\Longrightarrow \frac{%
k^{rn_{1}+n_{2}-(r+1)}-k^{n_{2}-1}}{n}\in {\Bbb N},  \label{43}
\end{equation}
from which

\[
\]
\begin{equation}
\frac{k^{rn_{1}+n_{2}-(r+1)}-1}{n}\in {\Bbb N}  \label{44}
\end{equation}
by using (24). The property (39) emerges just by iteration.

\smallskip

\smallskip

${\bf d)}$ It is easy to show that the property (39) is satisfied also for $%
\forall r,s\in Z$, provided that

\begin{equation}
rn_{1}+sn_{2}-(r+s)>0.  \label{47}
\end{equation}

{\bf e) }The property (39) can be further generalized to give

\smallskip

\begin{equation}
\frac{k^{rn_{1}^{q}+sn_{2}^{p}-(r+s)}-1}{n}\in {\Bbb N},\text{ }\forall
r,s\in {\Bbb N},\forall q,p\in {\Bbb N},  \label{53}
\end{equation}

Finally, this property can be extended in its turn assuming that $r,s\in Z$
and $rn_{1}^{q}+sn_{2}^{p}-(r+s)>0.$

\smallskip

\subsection{Pseudoprimes product of three primes}

\smallskip

$\bullet $ Let $n=n_{1}n_{2}n_{3}$ be a pseudoprime of basis $k\in {\Bbb N}%
:k>1,$ where $n_{1},n_{2},n_{3}$ are primes, and $k$ is prime relatively to $%
n_{1},n_{2},n_{3}.$ Then (see (15))

\begin{equation}
k^{n}-k=\Pi _{n}+\Pi _{n_{1}n_{2}}+\Pi _{n_{1}n_{3}}+\Pi _{n_{2}n_{3}}+\Pi
_{n_{1}}+\Pi _{n_{2}}+\Pi _{n_{3}}.  \label{54}
\end{equation}

\smallskip Furthermore, since $n$ is a pseudoprime, i.e.

\[
\frac{k^{n}-k}{n}\in {\Bbb N}, 
\]

and $\frac{\Pi _{n}}{n}\in {\Bbb N},$ where $\frac{\Pi _{n}}{n}$ is the
number of (independent) orbits of prime period $n,$ we have

\smallskip

\begin{equation}
\frac{\Pi _{n_{1}n_{2}}+\Pi _{n_{1}n_{3}}+\Pi _{n_{2}n_{3}}+\Pi _{n_{1}}+\Pi
_{n_{2}}+\Pi _{n_{3}}}{n}\in {\Bbb N}.  \label{56}
\end{equation}

\smallskip

$\bullet $ We notice that

\begin{equation}
\Pi _{n_{j}}=k^{n_{j}}-k,\text{ }j=1,2,3,  \label{57}
\end{equation}
and

\begin{equation}
\Pi
_{n_{i}n_{j}}=(k^{_{n_{i}n_{j}}}-k)-(k^{n_{i}}-k)-(k^{n_{j}}-k)=k^{n_{i}n_{j}}-k^{n_{i}}-k^{n_{j}}+k,
\label{58}
\end{equation}
for $i\neq j=1,2,3.$ Then, the condition (44) is equivalent to

\begin{equation}
\frac{%
\begin{array}{c}
\\ 
k^{n_{1}n_{2}}+k^{n_{1}n_{3}}+k^{n_{2}n_{3}}-k^{n_{1}}-k^{n_{2}}-k^{n_{3}}
\end{array}
}{n}\in {\Bbb N},  \label{59}
\end{equation}
which obviously implies 
\begin{equation}
\frac{%
\begin{array}{c}
\\ 
k^{n_{1}n_{2}}+k^{n_{1}n_{3}}+k^{n_{2}n_{3}}-k^{n_{1}}-k^{n_{2}}-k^{n_{3}}
\end{array}
}{n_{i}}\in {\Bbb N}  \label{60}
\end{equation}
for $i=1,2,3.$

\smallskip

$\bullet $ Since $n_{1}$ is a prime, the following properties

\smallskip

\begin{equation}
\frac{k^{n_{1}n_{2}}-k^{n_{2}}}{n_{1}}=\frac{(k^{n_{2}})^{n_{1}}-k^{n_{2}}}{%
n_{1}}\in {\Bbb N},  \label{61}
\end{equation}

\smallskip

\begin{equation}
\frac{k^{n_{1}n_{3}}-k^{n_{3}}}{n_{1}}\in {\Bbb N},  \label{62}
\end{equation}

hold.

\[
\]
Then, keeping in mind (48) and (49) and assuming that $n_{2}n_{3}>n_{1},$
from (47) we find

\begin{equation}
\frac{k^{n_{2}n_{3}}-k^{n_{1}}}{n_{1}}\in {\Bbb N}.  \label{63}
\end{equation}
In a similar way, assuming that $n_{1}n_{2}>n_{3}$ and $n_{1}n_{3}>n_{2},$
we obtain

\smallskip

\begin{equation}
\frac{k^{n_{1}n_{2}}-k^{n_{3}}}{n_{3}}\in {\Bbb N},  \label{64}
\end{equation}
and

\begin{equation}
\frac{k^{n_{1}n_{3}}-k^{n_{2}}}{n_{2}}\in {\Bbb N},  \label{65}
\end{equation}
respectively.

\smallskip

{\bf Remark}

\smallskip

The relations (52), (53) and (54) can be written as

\begin{equation}
\frac{k^{n_{2}n_{3}-n_{1}}-1}{n_{1}}\in {\Bbb N},  \label{66}
\end{equation}
\begin{equation}
\frac{k^{n_{1}n_{2}-n_{3}}-1}{n_{3}}\in {\Bbb N},
\end{equation}
\[
\]
\begin{equation}
\frac{k^{n_{1}n_{3}-n_{2}}-1}{n_{2}}\in {\Bbb N}.  \label{68}
\end{equation}

We note that if the order relations $n_{2}n_{3}>n_{1},$ and so on, turn out
to be not (fully or partially) satisfied, then (55), (56) and (57) must be
expressed as

\smallskip 
\begin{eqnarray}
\frac{k^{\mid n_{2}n_{3}-n_{1}\mid }-1}{n_{1}} &\in &{\Bbb N}, \\
&&  \nonumber
\end{eqnarray}
and so on.

Moreover, it is easy to show that, in general, we have

\begin{equation}
\frac{k_{{}}^{j\mid n_{2}n_{3}-n_{1}^{m}\mid }-1^{{}}}{n_{1}}\in {\Bbb N}%
,\forall m,j\in {\Bbb N},  \label{75}
\end{equation}

\begin{eqnarray*}
&& \\
&&
\end{eqnarray*}
\begin{equation}
\frac{k_{{}}^{j\mid n_{1}n_{2}-n_{3}^{m}\mid }-1^{{}}}{n_{3}}\in {\Bbb N}%
,\forall m,j\in {\Bbb N},  \label{76}
\end{equation}
and

\smallskip

\begin{equation}
\frac{k_{{}}^{j\mid n_{1}n_{3}-n_{2}^{m}\mid }-1^{{}}}{n_{2}}\in {\Bbb N}%
,\forall m,j\in {\Bbb N}.  \label{77}
\end{equation}

\section{\protect\smallskip Conclusions}

\smallskip We have exploited the principal features of the circle map to
obtain some properties of primes and pseudoprimes. On the one side, we have
reproduced a Fermat's Theorem on primes. A crucial role has been played by
the property $\Pi _{n}\equiv 0$ $(\func{mod}n),$ where $\Pi _{n}$ denotes
the number of periodic points of prime period $n$ in the circle map.
Moreover, we have proven an interesting property on pseudoprimes (Theorem
1). On the other side, we have considered pseudoprimes factorized into
products of two and three primes. In both cases, the circle map reveals as a
natural framework to study primes and pseudoprimes. The main results are
represented by Theorem 2 and by the properties (43) and (59)-(61).
Furthermore, we remark that the property (36) can be interpreted as a
consequence of Euler's Theorem.

Although we have limited ourselves to consider pseudoprimes which are
products of two and three primes, respectively, the formulas expressing the
properties achieved are simple enough to retain reliable an extension to
higher factorizations. We think also that the investigation carried out in
this note might be pursued by exploring the possible connection between
primes, pseudoprimes and other kinds of maps, as for example the logistic
map.

\smallskip

\smallskip

\smallskip

\smallskip

$^{1}$This work is supported in part by PRIN 97 ''SINTESI''

$^{2}$W. Scharlau, H. Opolka, {\it From Fermat to Minkowski }(Springer, New
York, 1985).

$^{3}$D.M. Bressoud, {\it Factorization and Primality Testing} (Springer,
New York, 1989).

$^{4}$See R.L. Liboff and M. Wong, ''Quasi-Chaotic Property of the
Prime-Number Sequence'', Int. J. Theor. Phys.{\bf 37}, 3109-3117\ (1998),
and references therein.

$^{5}$R.L. Devaney,{\it \ Chaotic Dynamical Systems} (Addison-Wesley,
Redwood City, California, 1989).

\[
\]

\end{document}